\documentclass{amsart}
\usepackage{amssymb}
\usepackage{mathrsfs}
\usepackage{hyperref}
\usepackage{verbatim}
\usepackage{graphicx}
\usepackage{subfigure}
\usepackage{psfrag}
\usepackage{enumitem}

\newtheorem{theorem}{Theorem}[section]

\newtheorem{proposition}[theorem]{Proposition}

{\theoremstyle{definition}

\newtheorem{example}[theorem]{Example}
}

\newcommand{\R}{\mathbb{R}}
\newcommand{\C}{\mathbb{C}}

\begin{document}
\title[The Pinchuk example revisited]{The Pinchuk example revisited}

\author[F. Fernandes \MakeLowercase{and} Z. Jelonek]
{Filipe Fernandes$^\dagger$ \MakeLowercase{and} Zbigniew Jelonek$^{*}$  }

\address{$^\dagger$ Departamento de Matem\'{a}tica, Universidade Federal do Cear\'a, 60440-900 Fortaleza, Cear\'a, Brazil}
\email{filipematonb@gmail.com} 

\address[Zbigniew Jelonek]{ Instytut Matematyczny, Polska Akademia Nauk, \'Sniadeckich 8, 00-656 Warszawa, Poland \& Departamento de Matem\'atica, Universidade Federal do Cear\'a,
	      Rua Campus do Pici, s/n, Bloco 914, Pici, 60440-900, 
	      Fortaleza-CE, Brazil. \newline  
              E-mail: {\tt najelone@cyf-kr.edu.pl}
}



\date{\today}

\begin{abstract}
In this note we provide two special examples of non-injective polynomial maps from $\R^2$ to $\R^2$ with non-vanishing Jacobian: the first one is  surjective,   the second one has non-dense image.
\end{abstract}

\maketitle

\section{Introduction} The famous Jacobian Conjecture is the following:

\vspace{5mm}

{\bf Jacobian Conjecture.} {\it Let $F:\C^2\to \C^2$ be a polynomial mapping with a constant non-zero Jacobian. Then $F$ is a bijection.}

\vspace{5mm}

There was also a similar real version of this conjecture:

\vspace{5mm}
{\bf Strong Real Jacobian Conjecture.} {\it Let $F:\R^2\to \R^2$ be a polynomial mapping with a  non-vanishing Jacobian. Then $F$ is a bijection.}

\vspace{5mm}

However, over twenty years ago Pinchuk in  \cite{P} found an example of locally injective polynomial endomorphism of $\R^2$, which is not injective. Pinchuk example has degree two and it is not a surjective mapping, in fact it omits two points. We can ask whether the Strong Real Jacobian Conjecture is true if we additionally assume that $F$ is surjective (such a question was posed by R. Peretz). Or is it possible to construct a counterexample of Pinchuk type such that the image $F(\R^2)$ is not dense in $\R^2$?

\section{Construction of the examples}
The main idea is to modify a existing non-injective polynomial map with non-vanishing Jacobian determinant by composing with special mappings. Here we use the version of Pinchuk map presented in \cite{BF}, because it seems to be the simplest known so far, in the sense that its components have the lowest degree (however this construction still can be adapted for any of the examples known, see \cite{F,P}). 
Let
\begin{align}
\begin{split}
p&=4\,{x}^{6}{y}^{3}+12\,{x}^{5}{y}^{2}+12\,{x}^{4}y+4\,{x}^{3}{y}^{2}+4
\,{x}^{3}+5\,{x}^{2}y+x+y, \\
q&= - \left( 4\,{x}^{4}{y}^{2}+8\,{x}^{3}y+4\,{x}^{2}+2\,xy+1 \right)\cdot  \\
& \cdot \left( 4\,{x}^{6}{y}^{3}+12\,{x}^{5}{y}^{2}+12\,{x}^{4}y+4\,{x}^{3}{y}^{2}+4\,{x}^{3}+7\,{x}^{2}y+3\,x+y \right) ,\label{qF}
\end{split}
\end{align}
we obtain the following:
\begin{theorem}\label{tF}
Let $p$ and $q$ be as in $\eqref{qF}$ and let $F:\mathbb{R}^2\to \mathbb{R}^2$ be defined by $F=(p,q)$. Then $F$ is a non-injective polynomial map with non-vanishing Jacobian and $(\pm1,0)$ are the only points with no inverse image.
\end{theorem}
\begin{proof}
See \cite[Theorem 3.3 and Corolary 4.6]{BF}.
\end{proof}

In the following proposition we present a surjective example.

\begin{proposition}\label{J}
Let $F$ be as in Theorem \ref{tF} and $\phi:\mathbb{C}\mapsto z^3-3z \in \mathbb{C}$, where we treat $\mathbb{C}$ as $\mathbb{R}^2$. Then $f=\phi \circ F$ is a non-injective, surjective polynomial map with non-vanishing Jacobian determinant.
\end{proposition}
\begin{proof}
The map $\phi$ is singular only for $z=\pm 1$. Since $(\pm1,0) \notin Im(F)$, $f$ has non-vanishing Jacobian determinant. Also $f$ can omit only the points $\phi(\pm1)=\mp 2$. However $\phi(\pm 2)=\pm2$, therefore $f$ is surjective.
\end{proof}

Here we present an example with non-dense image.

\begin{proposition}
Let 
\begin{equation}
\psi=((xy^2+x+y)(x-y),(xy+1)^2+x^2).
\end{equation}
and $\tilde{F}=(p+1,q)$, where $p$ and $q$ are as in \eqref{qF}. Then $\tilde{f}=\psi \circ \tilde{F}$ is a non-injective polynomial map with non-vanishing Jacobian determinant such that $\tilde{f}(\mathbb{R}^2)\subset \mathbb{R} \times \mathbb{R}^+$.
\end{proposition}
\begin{proof}
We have $$det(D(\psi))=(2y+2xy^2+x)^2+x^2(2y^2+3)^2,$$ thus $(0,0)$ is the only singular point of $\psi$. From Theorem \ref{tF} we have $(0,0) \notin Im(\tilde{F})$. Hence the Jacobian determinant of $\tilde{f}$ is non-vanishing. Also $(xy+1)^2+x^2>0$ for every $(x,y)\in \mathbb{R}^2$. Therefore $Im(\psi)\subset \mathbb{R} \times \mathbb{R}^+$ is not dense in $\mathbb{R}^2$. Consequently also $Im(\tilde{f})\subset  \mathbb{R} \times \mathbb{R}^+$.
\end{proof}

We finish this note answering negatively the following question of professor Tiep Dinh:

\vspace{5mm}
{\bf Question.}{\it Let $f: \mathbb{R}^n \to \mathbb{R}^n$ be a semi-algebraic smooth mapping which is surjective and it has a constant non-zero Jacobian. Then is $f$ a diffeomorphism?}

\vspace{5mm}

\begin{example}[Compare with \cite{LAC}]\label{exmp}
Let $F$ be as in Theorem \ref{J} and from every map $X:\mathbb{R}^n\to \mathbb{R}^n$ denote by $J(X)$ the determinant of the Jacobian of $X$. Then $G=\left(F,\frac{z}{J(F)}\right)$ is a rational non-injective local diffeomorphism of $\mathbb{R}^3$ with constant Jacobian determinant.
\end{example}
\begin{proof}
$G$ is well defined for all $(x,y,z)\in \mathbb{R}^3$, since $J(F)>0$. As for J(G), we have
$$J(G)=J\left(F,\frac{z}{J(F)}\right)=J(F,1/J(F))+J(F,z)/J(F)=J(F,z)/J(F)=1.$$
It is easy to see that $G$ is surjective.
\end{proof}

\section{Acknowledgments}
The first named author was partially supported by the grant 160946/2022-0 of Instituto Nacional de Ci\^encia e Tecnologia de Matem\'atica.

\end{document}